\documentclass[12pt]{amsart}
\usepackage{amssymb}
\usepackage{t1enc}
\usepackage[latin2]{inputenc}
\usepackage{verbatim}
\usepackage{amsmath,amsfonts,amssymb,amsthm}
\usepackage[mathcal]{eucal}
\usepackage{enumerate}
\usepackage[centertags]{amsmath}
\usepackage{graphics}

\newtheorem{theorem}{Theorem}
\newtheorem{lemma}{Lemma}

\numberwithin{equation}{subsection}

\begin{document}
\author{G. Tephnadze}
\title[Strong convergence ]{Strong convergence of two--dimensional
Walsh-Fourier series}
\address{G. Tephnadze, Department of Mathematics, Faculty of Exact and
Natural Sciences, Tbilisi State University, Chavchavadze str. 1, Tbilisi
0128, Georgia}
\email{giorgitephnadze@gmail.com}
\date{}
\maketitle

\begin{abstract}
We prove that certain mean of the quadratical partial sums of the
two-dimensional Walsh-Fourier series are uniformly bounded operators from
the Hardy space $H_{p}$ to the space $L_{p}$ for $0<p<1.$
\end{abstract}



\section{INTRODUCTION}

It is known $\left[ 7,\text{ p. 125}\right] $ that the Walsh-Paley system is
not a Schauder basis in $L_{1}\left( G\right) $. Moreover, (see \cite{S-W-S}%
) there exists a function in the dyadic Hardy space $H_{1}\left( G\right) $,
the partial sums of which are not bounded in $L_{1}\left( G\right) .$
However, in Simon \cite{Si} the following strong convergence result was
obtained for all $f\in H_{1}:$%
\begin{equation*}
\underset{n\rightarrow \infty }{\lim }\frac{1}{\log n}\overset{n}{\underset{%
k=1}{\sum }}\frac{\left\Vert S_{k}f-f\right\Vert _{1}}{k}=0,
\end{equation*}%
where $S_{k}f$ denotes the $k-$th partial sum of the Walsh-Fourier series of
$f$ (For the trigonometric analogue see Smith \cite{sm}, for the Vilenkin
system see Gát \cite{gat1}).

Simon \cite{si1} proved that\textbf{\ }there is an absolute constant $c_{p},$
depends only $p,$ such that
\begin{equation}
\overset{\infty }{\underset{k=1}{\sum }}\frac{\left\Vert S_{k}f\right\Vert
_{p}^{p}}{k^{2-p}}\leq c_{p}\left\Vert f\right\Vert _{H_{p}}^{p},
\label{1cc}
\end{equation}%
for all $f\in H_{p},$ where $0<p<1.$

The author \cite{tep5} proved that sequence $\left\{ 1/k^{2-p}\right\}
_{k=1}^{\infty }$ in inequality (\ref{1cc}) is important.

For the two-dimensional Walsh-Fourier series Weisz \cite{We} generalized the
result of Simon and proved that if $\alpha \geq 0$ and $f\in H_{p}\left(
G\times G\right) ,$ then

\begin{equation*}
\underset{n,m\geq 2}{\sup }\left( \frac{1}{\log n\log m}\right) ^{\left[ p%
\right] }\underset{2^{-\alpha }\leq k/l\leq 2^{\alpha },\text{ }\left(
k,l\right) \leq \left( n,m\right) }{\sum }\frac{\left\Vert
S_{k,l}f\right\Vert _{p}^{p}}{\left( kl\right) ^{2-p}}\leq c\left\Vert
f\right\Vert _{H_{p}}^{p},
\end{equation*}%
where $0<p<1$ and $\left[ p\right] $ denotes the integer part of $p.$

Goginava and Gogoladze \cite{gg} proved that the following result is true:

\textbf{Theorem G}.Let $f\in H_{1}\left( G\times G\right) $. Then there
exists absolute constant $c$, such that
\begin{equation*}
\sum\limits_{n=1}^{\infty }\frac{\left\Vert S_{n,n}f\right\Vert _{1}}{n\log
^{2}n}\leq c\left\Vert f\right\Vert _{H_{1}}.
\end{equation*}

For two-dimensional trigonometric system analogical theorem was proved in
\cite{Go}.

Convergence of quadratical partial sums of two-dimensional Walsh-Fourier
series was investigated in details by Weisz \cite{Webook2}, Goginava \cite%
{go}, Gát, Goginava, Nagy \cite{GGN}, Gát, Goginava, Tkebuchava \cite{GGT}.

The main aim of this paper is to prove (see Theorem 1) that
\begin{equation}
\sum\limits_{n=1}^{\infty }\frac{\left\Vert S_{n,n}f\right\Vert _{p}^{p}}{%
n^{3-2p}}\leq c_{p}\left\Vert f\right\Vert _{H_{p}}^{p},  \label{th}
\end{equation}%
for all $f\in H_{p}\left( G\times G\right) ,$ where $0<p<1.$ We also proved
that sequence $\left\{ 1/n^{3-2p}\right\} _{n=1}^{\infty }$ in inequality %
\ref{th} is important (see Theorem 2).

\section{DEFINITIONS AND NOTATIONS}

Let $\mathbf{P}$ denote the set of positive integers, $\mathbf{N:=P\cup \{}0%
\mathbf{\}.}$ Denote by $Z_{2}$ the discrete cyclic group of order 2, that
is $Z_{2}=\{0,1\},$ where the group operation is the modulo 2 addition and
every subset is open. The Haar measure on $Z_{2}$ is given such that the
measure of a singleton is 1/2. Let $G$ be the complete direct product of the
countable infinite copies of the compact groups $Z_{2}.$ The elements of $G$
are of the form $x=\left( x_{0},x_{1},...,x_{k},...\right) $ with $x_{k}\in
\{0,1\}\left( k\in \mathbf{N}\right) .$ The group operation on $G$ is the
coordinate-wise addition, the measure (denote\thinspace $\,$by$\,\,\mu $)
and the topology are the product measure and topology. The compact Abelian
group $G$ is called the Walsh group. A base for the neighborhoods of $G$ can
be given in the following way:
\begin{eqnarray*}
I_{0}\left( x\right) &:&=G,\,\,\,I_{n}\left( x\right) :=\,I_{n}\left(
x_{0},...,x_{n-1}\right) \\
&:&=\left\{ y\in G:\,y=\left( x_{0},...,x_{n-1},y_{n},y_{n+1},...\right)
\right\} ,
\end{eqnarray*}
\begin{equation*}
\,\left( x\in G,n\in \mathbf{N}\right) .
\end{equation*}
These sets are called the dyadic intervals. Let $0=\left( 0:i\in \mathbf{N}%
\right) \in G$ denote the null element of $G,\,\,\,I_{n}:=I_{n}\left(
0\right) \,\left( n\in \mathbf{N}\right) .$ Set $e_{n}:=\left(
0,...,0,1,0,...\right) \in G$ the $n\,$th\thinspace coordinate of which is 1
and the rest are zeros $\left( n\in \mathbf{N}\right) .$ Let $\overline{I}%
_{n}:=G\backslash I_{n}.$

If $n\in \mathbf{N}$, then $n=\sum\limits_{i=0}^{\infty }n_{i}2^{i},$ where $%
n_{i}\in \{0,1\}\,\,\left( i\in \mathbf{N}\right) $, i. e. $n$ is expressed
in the number system of base 2. Denote $\left| n\right| :=\max \{j\in
\mathbf{N:}n_{j}\neq 0\}$, that is, $2^{\left| n\right| }\leq n<2^{\left|
n\right| +1}.$

It is easy to show that for every odd number $n_{0}=1$ and we can write $%
n=1+\sum_{i=1}^{\left\vert n\right\vert }n_{j}2^{i}$, where $n_{j}\in
\left\{ 0,1\right\} ,$ $~(j\in \mathbf{P})$.

For $k\in \mathbf{N}$ and $x\in G$ let as denote by
\begin{equation*}
r_{k}\left( x\right) :=\left( -1\right) ^{x_{k}}\,\,\,\,\,\,\left( x\in
G,k\in \mathbf{N}\right)
\end{equation*}%
the $k$-th Rademacher function.

The Walsh-Paley system is defined as the sequence of Walsh-Paley functions:
\begin{equation*}
w_{n}\left( x\right) :=\prod\limits_{k=0}^{\infty }\left( r_{k}\left(
x\right) \right) ^{n_{k}}=r_{\left\vert n\right\vert }\left( x\right) \left(
-1\right) ^{\sum\limits_{k=0}^{\left\vert n\right\vert
-1}n_{k}x_{k}}\,\,\,\,\,\,\left( x\in G,n\in \mathbf{P}\right) .
\end{equation*}

The Walsh-Dirichlet kernel is defined by
\begin{equation*}
D_{n}\left( x\right) =\sum\limits_{k=0}^{n-1}w_{k}\left( x\right) .
\end{equation*}

Recall that (see $\left[ 8,\text{ p. 7}\right] $)
\begin{equation}
D_{2^{n}}\left( x\right) =\left\{
\begin{array}{c}
2^{n},\text{ \qquad }x\in I_{n} \\
0,\,\,\,\text{\qquad }x\in \overline{I}_{n}%
\end{array}%
\right. ,  \label{dir1}
\end{equation}

Furthermore, the following representation holds for the $D_{n}`$s. Let $n\in
\mathbf{N}$ and $n=\sum\limits_{i=0}^{\infty }n_{i}2^{i}.$ Then
\begin{equation}
D_{n}\left( x\right) =w_{n}\left( x\right) \sum\limits_{j=0}^{\infty
}n_{j}w_{2^{j}}\left( x\right) D_{2^{j}}\left( x\right) .  \label{dir2}
\end{equation}

The rectangular partial sums of the 2-dimensional Walsh-Fourier series of
function $f\in L_{2}\left( G\times G\right) $ are defined as follows:

\begin{equation*}
S_{M,N}f\left( x,y\right) :=\sum\limits_{i=0}^{M-1}\sum\limits_{j=0}^{N-1}%
\widehat{f}\left( i,j\right) w_{i}\left( x\right) w_{j}\left( y\right) ,
\end{equation*}
where the numbers
\begin{equation*}
\widehat{f}\left( i,j\right) =\int\limits_{G\times G}f\left( x,y\right)
w_{i}\left( x\right) w_{j}\left( y\right) d\mu \left( x,y\right)
\end{equation*}
is said to be the $\left( i,j\right) -$th Walsh-Fourier coefficient of the
function \thinspace $f.$

Denote
\begin{equation*}
S_{M}^{\left( 1\right) }f\left( x,y\right) :=\int\limits_{G}f\left(
s,y\right) D_{M}\left( x+s\right) d\mu \left( s\right)
\end{equation*}
and
\begin{equation*}
S_{N}^{\left( 2\right) }f\left( x,y\right) :=\int\limits_{G}f\left(
x,t\right) D_{N}\left( y+t\right) d\mu \left( t\right) .
\end{equation*}

The norm (or quasinorm) of the space $L_{p}(G\times G)$ is defined by \qquad

\begin{equation*}
\left\| f\right\| _{p}:=\left( \int_{G\times G}\left| f\right| ^{p}d\mu
\right) ^{1/p},\qquad \left( 0<p<\infty \right) .
\end{equation*}

The space $weak-L_{p}\left( G\times G\right) $ consists of all measurable
functions $f$ for which

\begin{equation*}
\left\| f\right\| _{weak-L_{p}(G\times G)}:=\underset{\lambda >0}{\sup }%
\lambda \mu \left( f>\lambda \right) ^{1/p}<+\infty .
\end{equation*}

The $\sigma -$algebra generated by the dyadic 2-dimensional $I_{n}\left(
x\right) \times I_{n}\left( y\right) $ square of measure $2^{-n}\times
2^{-n} $ will be denoted by $\digamma _{n,n}\left( n\in \mathbf{N}\right) .$
Denote by $f=\left( f_{n,n}\text{ }n\in \mathbf{N}\right) $ one-parameter
martingale with respect to $\digamma _{n,n}\left( n\in \mathbf{N}\right) .$%
(for details see e.g. \cite{Webook1}).

The expectation operator and the conditional expectation operator relative
to the $\digamma _{n,n}\left( n\in \mathbf{N}\right) $ are denoted by $E$
and $E_{n,n}$ , respectively.

The maximal function of a martingale $f$ is defined by

\begin{equation*}
f^{\ast }=\sup_{n\in \mathbf{N}}\left\vert f_{n,n}\right\vert .
\end{equation*}

Let $f\in L_{1}\left( G\times G\right) $. Then the dyadic maximal function
is given by
\begin{equation*}
f^{*}\left( x,y\right) =\sup\limits_{n\in \mathbf{N}}\frac{1}{\mu \left(
I_{n}(x)\times I_{n}(y)\right) }\left| \int\limits_{I_{n}(x)\times
I_{n}(y)}f\left( s,t\right) d\mu \left( s,t\right) \right| ,\,\,
\end{equation*}
\begin{equation*}
\left( x,y\right) \in G\times G.
\end{equation*}

The dyadic Hardy space $H_{p}(G\times G)$ $\left( 0<p<\infty \right) $
consists of all functions for which

\begin{equation*}
\left\| f\right\| _{H_{p}}:=\left\| f^{*}\right\| _{p}<\infty .
\end{equation*}

If $f\in L_{1}\left( G\times G\right) ,$ then it is easy to show that the
sequence $\left( S_{2^{n},2^{n}}\left( f\right) :n\in \mathbf{N}\right) $ is
a martingale. If $f=\left( f_{n,n},n\in \mathbf{N}\right) $ is a martingale,
then the Walsh-Fourier coefficients must be defined in a slightly different
manner: $\qquad \qquad $
\begin{equation*}
\widehat{f}\left( i,j\right) :=\lim_{k\rightarrow \infty
}\int_{G}f_{k,k}\left( x,y\right) w_{i}\left( x\right) w_{j}\left( y\right)
d\mu \left( x,y\right) .
\end{equation*}

It is known \cite{tep4} that that Fourier coefficients of $f\in H_{p}\left(
G\times G\right) $ are not bounded when $0<p<1.$

The Walsh-Fourier coefficients of $f\in L_{1}\left( G\times G\right) $ are
the same as those of the martingale $\left( S_{2^{n},2^{n}}f:n\in \mathbf{N}%
\right) $ obtained from $f$ .

A bounded measurable function $a$ is a p-atom, if there exists a
dyadic\thinspace 2-dimensional cube $I\times I\mathbf{,}$ such that
\begin{equation*}
\left\{
\begin{array}{l}
a)\qquad \int_{I\times I}ad\mu =0, \\
b)\ \qquad \left\| a\right\| _{\infty }\leq \mu (I\times I)^{-1/p}, \\
c)\qquad \text{supp}\left( a\right) \subset I\times I.\qquad%
\end{array}
\right.
\end{equation*}

\section{FORMULATION OF MAIN RESULTS}

\begin{theorem}
Let $0<p<1$ and $f\in H_{p}\left( G\times G\right) $. Then
\begin{equation*}
\sum\limits_{n=1}^{\infty }\frac{\left\| S_{n,n}f\right\| _{p}^{p}}{n^{3-2p}}%
\leq c_{p}\left\| f\right\| _{H_{p}}^{p}.
\end{equation*}
\end{theorem}

\begin{theorem}
Let $0<p<1$ and $\Phi :\mathbf{N}\rightarrow [1,$ $\infty )$ is any
nondecreasing function, satisfying the condition $\underset{n\rightarrow
\infty }{\lim }\Phi \left( n\right) =+\infty .$ Then there exists a
martingale $f\in H_{p}\left( G\times G\right) $ such that
\end{theorem}

\begin{equation*}
\underset{n=1}{\overset{\infty }{\sum }}\frac{\left\| S_{n,n}f\right\|
_{weak-L_{p}}^{p}\Phi \left( n\right) }{n^{3-2p}}=\infty .
\end{equation*}

\section{AUXILIARY PROPOSITIONS}

\begin{lemma}
\cite{Webook1} A martingale $f\in L_{p}\left( G\times G\right) $ is in $%
H_{p}\left( G\times G\right) \left( 0<p\leq 1\right) $ if and only if there
exist a sequence $\left( a_{k},k\in \mathbf{N}\right) $ of p-atoms and a
sequence $\left( \mu _{k},k\in \mathbf{N}\right) $ of a real numbers such
that
\begin{equation}
\qquad \sum_{k=0}^{\infty }\mu _{k}E_{n,n}a_{k}=f_{n,n}  \label{a1}
\end{equation}
and
\end{lemma}

\begin{equation*}
\qquad \sum_{k=0}^{\infty }\left| \mu _{k}\right| ^{p}<\infty ,
\end{equation*}
Moreover, $\left\| f\right\| _{H_{p}}\backsim \inf \left( \sum_{k=0}^{\infty
}\left| \mu _{k}\right| ^{p}\right) ^{1/p}$, where the infimum is taken over
all decomposition of $f$ of the form (\ref{a1}).

\section{PROOF OF THE THEOREMS}

\textbf{Proof of Theorem 1.} If we apply Lemma 1 we only have to prove that
\begin{equation}
\sum\limits_{n=1}^{\infty }\frac{\left\| S_{n,n}a\right\| _{p}^{p}}{n^{3-2p}}%
\leq c_{p}<\infty ,  \label{main}
\end{equation}
for every $p$ atom $a$.

Let $a$ be an arbitrary $p$-atom with support $I_{N}\left( z^{\prime
}\right) \times I_{N}\left( z^{\prime \prime }\right) $ and $\mu \left(
I_{N}\right) =\mu \left( I_{N}\right) =2^{-N}$. We can suppose that $%
z^{\prime }=z^{\prime \prime }=0.$

Let $\left( x,y\right) \in \overline{I}_{N}\times \overline{I}_{N}$. In this
case $D_{2^{i}}\left( x+s\right) 1_{I_{N}}\left( s\right) =0$ and $%
D_{2^{i}}\left( y+t\right) 1_{I_{N}}\left( t\right) =0$ for $i\geq N$.
Recall that $w_{2^{j}}\left( x+t\right) =w_{2^{j}}\left( x\right) $ for $%
t\in I_{N}$ and $j<N$. Consequently, from (\ref{dir2}) we obtain%
\begin{eqnarray*}
&&S_{n,n}a\left( x,y\right) \\
&=&\int\limits_{G\times G}a\left( s,t\right) D_{n}\left( x+s\right)
D_{n}\left( y+t\right) d\mu \left( s,t\right) \\
&=&\int\limits_{I_{N}\times I_{N}}a\left( s,t\right) D_{n}\left( x+s\right)
D_{n}\left( y+t\right) d\mu \left( s,t\right) \\
&=&\int\limits_{I_{N}\times I_{N}}a\left( s,t\right) w_{n}\left(
x+s+y+t\right) \sum\limits_{i=0}^{N-1}n_{i}w_{2^{i}}\left( x+s\right)
D_{2^{i}}\left( x+s\right) \\
&&\times \sum\limits_{j=0}^{N-1}n_{j}w_{2^{j}}\left( y+t\right)
D_{2^{j}}\left( y+t\right) d\mu \left( s,t\right) \\
&=&w_{n}\left( x\right) \sum\limits_{i=0}^{N-1}n_{i}w_{2^{i}}\left( x\right)
D_{2^{i}}\left( x\right) w_{n}\left( y\right)
\sum\limits_{j=0}^{N-1}n_{j}w_{2^{j}}\left( y\right) D_{2^{j}}\left( y\right)
\\
&&\times \int\limits_{I_{N}\times I_{N}}a\left( s,t\right) w_{n}\left(
s+t\right) d\mu \left( s,t\right)
\end{eqnarray*}
\begin{eqnarray*}
&=&w_{n}\left( x+y\right) \sum\limits_{i=0}^{N-1}n_{i}w_{2^{i}}\left(
x\right) D_{2^{i}}\left( x\right)
\sum\limits_{j=0}^{N-1}n_{j}w_{2^{j}}\left( y\right) D_{2^{j}}\left( y\right)
\\
&&\times \int\limits_{I_{N}}\left( \int\limits_{I_{N}}a\left( t+\tau
,t\right) d\mu \left( t\right) \right) w_{n}\left( \tau \right) d\mu \left(
\tau \right) \\
&=&w_{n}\left( x+y\right) \sum\limits_{i=0}^{N-1}n_{i}w_{2^{i}}\left(
x\right) D_{2^{i}}\left( x\right)
\sum\limits_{j=0}^{N-1}n_{j}w_{2^{j}}\left( y\right) D_{2^{j}}\left(
y\right) \int\limits_{I_{N}}\Phi \left( \tau \right) w_{n}\left( \tau
\right) d\mu \left( \tau \right) \\
&=&w_{n}\left( x+y\right) \sum\limits_{i=0}^{N-1}n_{i}w_{2^{i}}\left(
x\right) D_{2^{i}}\left( x\right)
\sum\limits_{j=0}^{N-1}n_{j}w_{2^{j}}\left( y\right) D_{2^{j}}\left(
y\right) \widehat{\Phi }\left( n\right) ,
\end{eqnarray*}

where
\begin{equation*}
\Phi \left( \tau \right) =\int\limits_{I_{N}}a\left( t+\tau ,t\right) d\mu
\left( t\right) .
\end{equation*}

Let $x\in I_{s}\backslash I_{s+1}.$ Using (\ref{dir1}) we get
\begin{equation*}
\overset{N-1}{\underset{i=0}{\sum }}D_{2^{i}}\left( x\right) \leq c2^{s}.
\end{equation*}

Since
\begin{equation*}
\overset{-}{I_{N}}=\overset{N-1}{\underset{s=0}{\bigcup }}I_{s}\backslash
I_{s+1}
\end{equation*}
we obtain
\begin{eqnarray}
&&\int_{\overline{I}_{N}}\left( \overset{N-1}{\underset{i=0}{\sum }}%
D_{2^{i}}\left( x\right) \right) ^{p}d\mu \left( x\right)  \label{dir3} \\
&\leq &c_{p}\underset{s=0}{\overset{N-1}{\sum }}\int_{I_{s}\backslash
I_{s+1}}2^{ps}d\mu \left( x\right)  \notag \\
&\leq &c_{p}\underset{s=0}{\overset{\infty }{\sum }}2^{\left( p-1\right) s}
\notag \\
&<&c_{p}<\infty ,\text{ \qquad }0<p<1.  \notag
\end{eqnarray}
applying (\ref{dir3}) we can write
\begin{eqnarray*}
&&\sum\limits_{n=1}^{\infty }\frac{1}{n^{3-2p}}\int\limits_{\overline{I}%
_{N}\times \overline{I}_{N}}\left| S_{n,n}a\left( x,y\right) \right|
^{p}d\mu \left( x,y\right) \\
&\leq &\sum\limits_{n=1}^{\infty }\frac{\left| \widehat{\Phi }\left(
n\right) \right| ^{p}}{n^{3-2p}}\left( \int\limits_{\overline{I}_{N}}\left(
\sum\limits_{i=0}^{N-1}D_{2^{i}}\left( x\right) \right) ^{p}d\mu \left(
x\right) \right) ^{2} \\
&\leq &c_{p}\sum\limits_{n=1}^{\infty }\frac{\left| \widehat{\Phi }\left(
n\right) \right| ^{p}}{n^{3-2p}}.
\end{eqnarray*}

Let $n<2^{N}$. Since $w_{n}\left( \tau \right) =1,$ for $\tau \in I_{N}$ we
have
\begin{eqnarray*}
\widehat{\Phi }\left( n\right) &=&\int\limits_{I_{N}}\Phi \left( \tau
\right) w_{n}\left( \tau \right) d\mu \left( \tau \right) \\
&=&\int\limits_{I_{N}}\left( \int\limits_{I_{N}}a\left( t+\tau ,t\right)
d\mu \left( t\right) \right) w_{n}\left( \tau \right) d\mu \left( \tau
\right) \\
&=&\int\limits_{I_{N}\times I_{N}}a\left( s,t\right) d\mu \left( s,t\right)
=0.
\end{eqnarray*}
Hence, we can suppose that $n\geq 2^{N}.$ By Hölder inequality we obtain
\begin{eqnarray}
&&\sum\limits_{n=1}^{\infty }\frac{\left| \widehat{\Phi }\left( n\right)
\right| ^{p}}{n^{3-2p}}  \label{1p} \\
&\leq &\left( \sum\limits_{n=2^{N}}^{\infty }\left| \widehat{\Phi }\left(
n\right) \right| ^{2}\right) ^{p/2}\left( \sum\limits_{n=2^{N}}^{\infty }%
\frac{1}{n^{\left( 3-2p\right) \cdot \left( 2/\left( 2-p\right) \right) }}%
\right) ^{\left( 2-p\right) /2}  \notag \\
&\leq &\left( \frac{1}{2^{N\left( 2\left( 3-2p\right) /\left( 2-p\right)
-1\right) }}\right) ^{\left( 2-p\right) /2}\left( \int\limits_{G}\left| \Phi
\left( \tau \right) \right| ^{2}d\mu \left( \tau \right) \right) ^{p/2}
\notag \\
&\leq &\frac{c_{p}}{2^{N\left( 4-3p\right) /2}}\left(
\int\limits_{I_{N}}\left| \int\limits_{I_{N}}a\left( t+\tau ,t\right) d\mu
\left( t\right) \right| ^{2}d\mu \left( \tau \right) \right) ^{p/2}  \notag
\\
&\leq &\frac{c_{p}}{2^{N\left( 4-3p\right) /2}}\left\| a\right\| _{\infty
}^{p}\frac{1}{2^{Np/2}}\frac{1}{2^{Np}}  \notag \\
&\leq &\frac{c_{p}}{2^{N\left( 4-3p\right) /2}}2^{2N}\frac{1}{2^{3pN/2}}%
<c_{p}<\infty .  \notag
\end{eqnarray}

Let $\left( x,y\right) \in \overline{I}_{N}\times I_{N}$. Then we have
\begin{eqnarray*}
&&S_{n,n}a\left( x,y\right) \\
&=&w_{n}\left( x\right) \sum\limits_{j=0}^{N-1}n_{j}w_{2^{j}}\left( x\right)
D_{2^{j}}\left( x\right) \\
&&\times \int\limits_{G\times G}a\left( s,t\right) w_{n}\left( s\right)
D_{n}\left( y+t\right) d\mu \left( s,t\right) \\
&=&w_{n}\left( x\right) \sum\limits_{j=0}^{N-1}n_{j}w_{2^{j}}\left( x\right)
D_{2^{j}}\left( x\right) \int\limits_{G}S_{n}^{\left( 2\right) }a\left(
s,y\right) w_{n}\left( s\right) d\mu \left( s\right) \\
&=&w_{n}\left( x\right) \sum\limits_{j=0}^{N-1}n_{j}w_{2^{j}}\left( x\right)
D_{2^{j}}\left( x\right) \widehat{S}_{n}^{\left( 2\right) }a\left(
n,y\right) .
\end{eqnarray*}
Using (\ref{dir3}) we get
\begin{eqnarray*}
&&\sum\limits_{n=1}^{\infty }\frac{1}{n^{3-2p}}\int\limits_{\overline{I}%
_{N}\times I_{N}}\left| S_{n,n}a\left( x,y\right) \right| ^{p}d\mu \left(
x,y\right) \\
&\leq &\sum\limits_{n=1}^{\infty }\frac{1}{n^{3-2p}}\int\limits_{\overline{I}%
_{N}\times I_{N}}\left( \sum\limits_{j=0}^{N-1}D_{2^{j}}\left( x\right)
\left| \widehat{S}_{n}^{\left( 2\right) }a\left( n,y\right) \right| \right)
^{p}d\mu \left( x,y\right) \\
&\leq &\sum\limits_{n=1}^{\infty }\frac{1}{n^{3-2p}}\int_{\overline{I}%
_{N}}\left( \overset{N-1}{\underset{i=0}{\sum }}D_{2^{i}}\left( x\right)
\right) ^{p}d\mu \left( x\right) \cdot \int\limits_{I_{N}}\left| \widehat{S}%
_{n}^{\left( 2\right) }a\left( n,y\right) \right| ^{p}d\mu \left( y\right) \\
&\leq &\sum\limits_{n=1}^{\infty }\frac{1}{n^{3-2p}}\int\limits_{I_{N}}%
\left| \widehat{S}_{n}^{\left( 2\right) }a\left( n,y\right) \right| ^{p}d\mu
\left( y\right) .
\end{eqnarray*}

Let $n<2^{N}$. Then by the definition of the atom we have
\begin{eqnarray*}
\widehat{S}_{n}^{\left( 2\right) }a\left( n,y\right)
&=&\int\limits_{G}\left( \int\limits_{G}a\left( s,t\right) D_{n}\left(
y+t\right) d\mu \left( t\right) \right) w_{n}\left( s\right) d\mu \left(
s\right) \\
&=&D_{n}\left( y\right) \int\limits_{I_{N}\times I_{N}}a\left( s,t\right)
d\mu \left( s,t\right) =0.
\end{eqnarray*}%
Therefore, we can suppose that $n\geq 2^{N}$. Hence
\begin{eqnarray*}
&&\sum\limits_{n=1}^{\infty }\frac{1}{n^{3-2p}}\int\limits_{\overline{I}%
_{N}\times I_{N}}\left\vert S_{n,n}a\left( x,y\right) \right\vert ^{p}d\mu
\left( x,y\right) \\
&\leq &\sum\limits_{n=2^{N}}^{\infty }\frac{1}{n^{3-2p}}\int\limits_{I_{N}}%
\left\vert \widehat{S}_{n}^{\left( 2\right) }a\left( n,y\right) \right\vert
^{p}d\mu \left( y\right)
\end{eqnarray*}%
Since
\begin{equation*}
\left\Vert S_{n}^{\left( 2\right) }a\left( n,y\right) \right\Vert _{2}\leq
c\left\Vert a\right\Vert _{2}
\end{equation*}%
from Hölder inequality we can write%
\begin{eqnarray*}
&&\int\limits_{I_{N}}\left\vert \widehat{S}_{n}^{\left( 2\right) }a\left(
n,y\right) \right\vert ^{p}d\mu \left( y\right) \\
&\leq &\frac{c_{p}}{2^{N\left( 1-p\right) }}\left(
\int\limits_{I_{N}}\left\vert \widehat{S}_{n}^{\left( 2\right) }a\left(
n,y\right) \right\vert d\mu \left( y\right) \right) ^{p} \\
&=&\frac{c_{p}}{2^{N\left( 1-p\right) }}\left( \int\limits_{I_{N}}\left\vert
\int\limits_{I_{N}}S_{n}^{\left( 2\right) }a\left( s,y\right) w_{n}\left(
s\right) d\mu \left( s\right) \right\vert d\mu \left( y\right) \right) ^{p}
\end{eqnarray*}
\begin{eqnarray*}
&=&\frac{c_{p}}{2^{N\left( 1-p\right) }}\left( \int\limits_{I_{N}}\left\vert
\int\limits_{I_{N}}\left( \int\limits_{I_{N}}a\left( s,t\right) D_{n}\left(
y+t\right) d\mu \left( t\right) \right) w_{n}\left( s\right) d\mu \left(
s\right) \right\vert d\mu \left( y\right) \right) ^{p} \\
&\leq &\frac{c_{p}}{2^{N\left( 1-p\right) }}\left( \int\limits_{I_{N}}\left(
\int\limits_{I_{N}}\left\vert \int\limits_{I_{N}}a\left( s,t\right)
D_{n}\left( y+t\right) d\mu \left( t\right) \right\vert d\mu \left( y\right)
\right) d\mu \left( s\right) \right) ^{p} \\
&\leq &\frac{c_{p}}{2^{N\left( 1-p\right) }}\left( \frac{1}{2^{N/2}}%
\int\limits_{I_{N}}\left( \int\limits_{I_{N}}\left\vert
\int\limits_{I_{N}}a\left( s,t\right) D_{n}\left( y+t\right) d\mu \left(
t\right) \right\vert ^{2}d\mu \left( y\right) \right) ^{1/2}d\mu \left(
s\right) \right) ^{p} \\
&\leq &\frac{c_{p}}{2^{N\left( 1-p\right) }}\left( \frac{1}{2^{N/2}}%
\int\limits_{I_{N}}\left( \int\limits_{I_{N}}\left\vert a\left( s,t\right)
\right\vert ^{2}d\mu \left( t\right) \right) ^{1/2}d\mu \left( s\right)
\right) ^{p} \\
&\leq &\frac{c_{p}}{2^{N\left( 1-p\right) }}\left( \frac{\left\Vert
a\right\Vert _{\infty }}{2^{N/2}}\frac{1}{2^{N}}\frac{1}{2^{N/2}}\right)
^{p}\leq \frac{c_{p}}{2^{N\left( 1-p\right) }}\left( \frac{2^{2N/p}}{2^{2N}}%
\right) ^{p}\leq c_{p}2^{N\left( 1-p\right) }.
\end{eqnarray*}%
Consequently,
\begin{eqnarray}
&&\sum\limits_{n=1}^{\infty }\frac{1}{n^{3-2p}}\int\limits_{\overline{I}%
_{N}\times I_{N}}\left\vert S_{n,n}a\left( x,y\right) \right\vert d\mu
\left( x,y\right)  \label{2p} \\
&\leq &c_{p}\sum\limits_{n=2^{N}}^{\infty }\frac{1}{n^{3-2p}}2^{N\left(
1-p\right) }\leq \frac{c_{p}}{2^{N\left( 1-p\right) }}\leq c_{p}<\infty .
\notag
\end{eqnarray}

Analogously, we can prove that
\begin{equation}
\sum\limits_{n=1}^{\infty }\frac{1}{n^{3-2p}}\int\limits_{I_{N}\times
\overline{I}_{N}}\left| S_{n,n}a\left( x,y\right) \right| ^{p}d\mu \left(
x,y\right) \leq c_{p}<\infty .  \label{3p}
\end{equation}

Let $\left( x,y\right) \in I_{N}\times I_{N}$. Then by the definition of the
atom we can write
\begin{eqnarray*}
&&\int\limits_{I_{N}\times I_{N}}\left| S_{n,n}a\left( x,y\right) \right|
^{p}d\mu \left( x,y\right) \\
&\leq &\frac{1}{2^{N\left( 2-p\right) }}\left( \int\limits_{I_{N}\times
I_{N}}\left| S_{n,n}a\left( x,y\right) \right| ^{2}d\mu \left( x,y\right)
\right) ^{p/2} \\
&\leq &\frac{1}{2^{N\left( 2-p\right) }}\left( \int\limits_{I_{N}\times
I_{N}}\left| a\left( x,y\right) \right| ^{2}d\mu \left( x,y\right) \right)
^{p/2} \\
&\leq &\frac{\left\| a\right\| _{\infty }^{p}}{2^{N\left( 2-p\right) }}\frac{%
1}{2^{Np}}\leq c_{p}\frac{1}{2^{N\left( 2-p\right) }}2^{2N}\frac{1}{2^{Np}}%
\leq c_{p}<\infty .
\end{eqnarray*}

It follows that
\begin{eqnarray}
&&\sum\limits_{n=1}^{\infty }\frac{1}{n^{3-2p}}\int\limits_{I_{N}\times
I_{N}}\left| S_{n,n}a\left( x,y\right) \right| d\mu \left( x,y\right)
\label{4p} \\
&\leq &c_{p}\sum\limits_{n=1}^{\infty }\frac{1}{n^{3-2p}}\leq c_{p}<\infty .
\notag
\end{eqnarray}

Combining (\ref{main}-\ref{4p}) we complete the proof of Theorem 1.

\textbf{Proof of Theorem 2. }Let $0<p<1$ and $\Phi \left( n\right) $ is any
nondecreasing, nonnegative function, satisfying condition
\begin{equation*}
\underset{n\rightarrow \infty }{\lim }\Phi \left( n\right) =\infty ,
\end{equation*}
For this function $\Phi \left( n\right) ,$ there exists an increasing
sequence of the positive integers $\left\{ \alpha _{k}:\text{ }k\geq
0\right\} $ such that:

\begin{equation*}
\alpha _{0}\geq 2
\end{equation*}
and

\begin{equation}
\sum_{k=0}^{\infty }\frac{1}{\Phi ^{p/4}\left( 2^{\alpha _{k}}\right) }%
<\infty .  \label{2}
\end{equation}

Let \qquad
\begin{equation*}
f_{A,A}\left( x,y\right) =\sum_{\left\{ k;\text{ }\alpha _{k}<A\right\}
}\lambda _{k}a_{k},
\end{equation*}
where
\begin{equation*}
\lambda _{k}=\frac{1}{\Phi ^{1/4}\left( 2^{\alpha _{k}}\right) }
\end{equation*}
and

\begin{equation*}
a_{k}\left( x,y\right) =2^{\alpha _{k}\left( 2/p-2\right) }\left(
D_{2^{\alpha _{k}+1}}\left( x\right) -D_{2^{\alpha _{k}}}\left( x\right)
\right) \left( D_{2^{\alpha _{k}+1}}\left( y\right) -D_{2^{\alpha
_{k}}}\left( y\right) \right) .
\end{equation*}

It is easy to show that the martingale $\,f=\left( f_{1,1},\text{ }%
f_{2,2},...,\text{ }f_{A,A},\text{ }...\right) \in H_{p}.$

Indeed, since

\begin{equation}
S_{2^{A}}a_{k}\left( x,y\right) =\left\{
\begin{array}{l}
a_{k}\left( x,y\right) \text{, \qquad }\alpha _{k}<A, \\
0\text{, \qquad }\alpha _{k}\geq A,%
\end{array}
\right.  \label{4}
\end{equation}

\begin{eqnarray*}
\text{supp}(a_{k}) &=&I_{\alpha _{k}}, \\
\int_{I_{\alpha _{k}}}a_{k}d\mu &=&0
\end{eqnarray*}
and

\begin{equation*}
\left\| a_{k}\right\| _{\infty }\leq 2^{\alpha _{k}\left( 2/p-2\right)
}2^{2\alpha _{k}}\leq 2^{2\alpha _{k}/p}=(\text{supp }a_{k})^{-1/p}
\end{equation*}
from Lemma 1 and (\ref{2}) we conclude that $f\in H_{p}.$

It is easy to show that

\begin{equation}
\widehat{f}(i,j)=\left\{
\begin{array}{l}
\frac{2^{\alpha _{k}\left( 2/p-2\right) }}{\Phi ^{1/4}\left( 2^{\alpha
_{k}}\right) },\,\, \\
\text{ if }\left( i,\text{\thinspace \thinspace }j\right) \in \left\{
2^{\alpha _{k}},...,\text{ ~}2^{\alpha _{k}+1}-1\right\} \times \left\{
2^{\alpha _{k}},...,\text{ ~}2^{\alpha _{k}+1}-1\right\} ,\text{ }k=0,1,2...
\\
0,\text{ \thinspace } \\
\text{\thinspace \thinspace if \thinspace }\left( \,i,\text{ \thinspace }%
j\right) \notin \bigcup\limits_{k=1}^{\infty }\left\{ 2^{\alpha _{k}},...,%
\text{ ~}2^{\alpha _{k}+1}-1\right\} \text{ }\times \left\{ 2^{\alpha
_{k}},...,\text{ ~}2^{\alpha _{k}+1}-1\right\} .%
\end{array}%
\right.  \label{5}
\end{equation}

Let\textbf{\ } $2^{\alpha _{k}}<n<2^{\alpha _{k}+1}$. From (\ref{5}) we have

\begin{eqnarray}
&&S_{n,n}f\left( x,y\right)  \label{13d} \\
&=&\sum_{i=0}^{2^{\alpha _{k-1}+1}-1}\sum_{j=0}^{2^{\alpha _{k-1}+1}-1}%
\widehat{f}(i,j)w_{i}\left( x\right) w_{j}\left( y\right)  \notag \\
&&+\sum_{i=2^{\alpha _{k}}}^{n-1}\sum_{j=2^{\alpha _{k}}}^{n-1}\widehat{f}%
(i,j)w_{i}\left( x\right) w_{j}\left( y\right)  \notag \\
&=&\sum_{\eta =0}^{k-1}\sum_{i=2^{\alpha _{\eta }}}^{2^{\alpha _{\eta
}+1}-1}\sum_{j=2^{\alpha _{\eta }}}^{2^{\alpha _{\eta }+1}-1}\widehat{f}%
(i,j)w_{i}\left( x\right) w_{j}\left( y\right)  \notag \\
&&+\sum_{i=2^{\alpha _{k}}}^{n-1}\sum_{j=2^{\alpha _{k}}}^{n-1}\widehat{f}%
(i,j)w_{i}\left( x\right) w_{j}\left( y\right)  \notag \\
&=&\sum_{\eta =0}^{k-1}\sum_{i=2^{\alpha _{\eta }}}^{2^{\alpha _{\eta
}+1}-1}\sum_{j=2^{\alpha _{\eta }}}^{2^{\alpha _{\eta }+1}-1}\frac{2^{\alpha
_{\eta }\left( 2/p-2\right) }}{\Phi ^{1/4}\left( 2^{\alpha _{\eta }}\right) }%
w_{i}\left( x\right) w_{j}\left( y\right)  \notag \\
&&+\sum_{i=2^{\alpha _{k}}}^{n-1}\sum_{j=2^{\alpha _{k}}}^{n-1}\frac{%
2^{\alpha _{k}\left( 2/p-2\right) }}{\Phi ^{1/4}\left( 2^{\alpha
_{k}}\right) }w_{i}\left( x\right) w_{j}\left( y\right)  \notag \\
&=&\sum_{\eta =0}^{k-1}\frac{2^{\alpha _{\eta }\left( 2/p-2\right) }}{\Phi
^{1/4}\left( 2^{\alpha _{\eta }}\right) }\left( D_{2^{\alpha _{\eta
}+1}}\left( x\right) -D_{2^{\alpha _{\eta }}}\left( x\right) \right) \left(
D_{2^{\alpha _{\eta }+1}}\left( y\right) -D_{2^{\alpha _{\eta }}}\left(
y\right) \right)  \notag \\
&&+\frac{2^{\alpha _{k}\left( 2/p-2\right) }}{\Phi ^{1/4}\left( 2^{\alpha
_{k}}\right) }\left( D_{_{n}}\left( x\right) -D_{2^{\alpha _{k}}}\left(
x\right) \right) \left( D_{n}\left( y\right) -D_{2^{\alpha _{k}}}\left(
y\right) \right)  \notag \\
&=&I+II.  \notag
\end{eqnarray}

Let $\left( x,y\right) \in \left( G\backslash I_{1}\right) \times \left(
G\backslash I_{1}\right) $ and $n\ $is odd number. Since $n-2^{\alpha _{k}}\
$is odd number too and
\begin{equation*}
D_{n+2^{\alpha _{k}}}\left( x\right) =D_{2^{\alpha _{k}}}\left( x\right)
+w_{2^{\alpha _{k}}}\left( x\right) D_{n}\left( x\right) ,\text{ when }%
\,\,n<2^{\alpha _{k}},
\end{equation*}%
from (\ref{dir1}) and (\ref{dir2}) we can write

\begin{eqnarray}
\left| II\right| &=&\frac{2^{\alpha _{k}\left( 2/p-2\right) }}{\Phi
^{1/4}\left( 2^{\alpha _{k}}\right) }\left| w_{2^{\alpha _{k}}}\left(
x\right) D_{n-2^{\alpha _{k}}}\left( x\right) w_{2^{\alpha _{k}}}\left(
y\right) D_{n-2^{\alpha _{k}}}\left( y\right) \right|  \label{13a} \\
&=&\frac{2^{\alpha _{k}\left( 2/p-2\right) }}{\Phi ^{1/4}\left( 2^{\alpha
_{k}}\right) }\left| w_{2^{\alpha _{k}}}\left( x\right) w_{n-2^{\alpha
_{k}}}\left( x\right) D_{_{1}}\left( x\right) w_{2^{\alpha _{k}}}\left(
y\right) w_{n-2^{\alpha _{k}}}\left( y\right) D_{_{1}}\left( y\right) \right|
\notag \\
&=&\frac{2^{\alpha _{k}\left( 2/p-2\right) }}{\Phi ^{1/4}\left( 2^{\alpha
_{k}}\right) }.  \notag
\end{eqnarray}
Applying (\ref{dir1}) and condition $\alpha _{n}\geq 2$ $\left( n\in \mathbf{%
N}\right) $ for $I$ we have

\begin{equation}
I=\sum_{\eta =0}^{k-1}\frac{2^{\alpha _{k}\left( 2/p-2\right) }}{\Phi
^{1/4}\left( 2^{\alpha _{\eta }}\right) }\left( D_{2^{\alpha _{\eta
}+1}}\left( x\right) -D_{2^{\alpha _{\eta }}}\left( x\right) \right) \left(
D_{2^{\alpha _{\eta }+1}}\left( y\right) -D_{2^{\alpha _{\eta }}}\left(
y\right) \right) =0.  \label{13b}
\end{equation}
Hence

\begin{eqnarray}
&&\left\| S_{n,n}f\left( x,y\right) \right\| _{weak-L_{p}}  \label{13} \\
&\geq &\frac{2^{\alpha _{k}\left( 2/p-2\right) }}{2\Phi ^{1/4}\left(
2^{\alpha _{k}}\right) }\left( \mu \left\{ \left( x,y\right) \in \left(
G\backslash I_{1}\right) \times \left( G\backslash I_{1}\right) :\left|
S_{n,n}f\left( x,y\right) \right| \geq \frac{2^{\alpha _{k}\left(
2/p-2\right) }}{2\Phi ^{1/4}\left( 2^{\alpha _{k}}\right) }\right\} \right)
^{1/p}  \notag \\
&\geq &\frac{2^{\alpha _{k}\left( 2/p-2\right) }}{2\Phi ^{1/4}\left(
2^{\alpha _{k}}\right) }\left| \left( G\backslash I_{1}\right) \times \left(
G\backslash I_{1}\right) \right| \geq \frac{c_{p}2^{\alpha _{k}\left(
2/p-2\right) }}{\Phi ^{1/4}\left( 2^{\alpha _{k}}\right) }.  \notag
\end{eqnarray}

Using (\ref{13}) we have

\begin{eqnarray}
&&\underset{n=1}{\overset{2^{\alpha _{k}+1}-1}{\sum }}\frac{\left\|
S_{n,n}f\right\| _{weak-L_{p}}^{p}\Phi \left( n\right) }{n^{3-2p}}
\label{14} \\
&\geq &\underset{n=2^{\alpha _{k}}+1}{\overset{2^{\alpha _{k}+1}-1}{\sum }}%
\frac{\left\| S_{n,n}f\right\| _{weak-L_{p}}^{p}\Phi \left( n\right) }{%
n^{3-2p}}  \notag \\
&\geq &c_{p}\Phi \left( 2^{\alpha _{k}}\right) \underset{n=2^{\alpha
_{k}-1}+1}{\overset{2^{\alpha _{k}}-1}{\sum }}\frac{\left\|
S_{2n+1,2n+1}f\right\| _{weak-L_{p}}^{p}}{\left( 2n+1\right) ^{3-2p}}  \notag
\\
&\geq &c_{p}\Phi \left( 2^{\alpha _{k}}\right) \frac{2^{2\alpha _{k}\left(
1-p\right) }}{\Phi ^{1/4}\left( 2^{\alpha _{k}}\right) }\underset{%
n=2^{\alpha _{k}-1}+1}{\overset{2^{\alpha _{k}}-1}{\sum }}\frac{1}{\left(
2n+1\right) ^{3-2p}}  \notag \\
&\geq &c_{p}\Phi ^{3/4}\left( 2^{\alpha _{k}}\right) \rightarrow \infty ,%
\text{ \qquad when }k\rightarrow \infty .  \notag
\end{eqnarray}

Combining (\ref{2}-\ref{14}) we complete the proof of Theorem 2.

\end{document}